\documentclass[a4paper]{amsart}
\usepackage{amsthm,amssymb}
\usepackage[T1]{fontenc}
\usepackage[utf8]{inputenc}
\usepackage{enumitem}
\usepackage{mathrsfs}
\usepackage{dsfont}
\usepackage{lmodern}
\usepackage{latexsym}
\usepackage{scalerel,stackengine}
\usepackage{nicefrac}
\usepackage{natbib}

\usepackage{caption}
\usepackage{graphicx}
\usepackage{hyperref}

\usepackage[dvipsnames]{xcolor}
\usepackage{tikz}
\usepackage{tikz-cd}
\usetikzlibrary{positioning,decorations.pathmorphing,intersections}

\tikzset{
myarrow/.style = {-stealth,ultra thick,shorten >=2pt,shorten <=2pt,cyan}
}

\tikzset{%
    symbol/.style={%
        draw=none,
        every to/.append style={%
            edge node={node [sloped, allow upside down, auto=false]{$#1$}}
            }
        }
    }

\numberwithin{equation}{section}

\newcommand\frB{\mathfrak{B}}

\DeclareMathSymbol{\ppartof}{\mathord}{AMSa}{64}
\DeclareMathSymbol{\partof}{\mathrel}{AMSa}{64}
\DeclareSymbolFont{symbolsC}{U}{txsyc}{m}{n}
\DeclareMathSymbol{\poverl}{\mathord}{symbolsC}{7}
\DeclareMathSymbol{\overl}{\mathrel}{symbolsC}{7}
\DeclareMathSymbol{\pext}{\mathord}{symbolsC}{78}
\DeclareMathSymbol{\npartof}{\mathrel}{symbolsC}{97}
\DeclareMathSymbol{\ningr}{\mathrel}{symbolsC}{64}
\DeclareMathSymbol{\nll}{\mathrel}{symbolsC}{51}

\newcommand\ext{\mathrel{\bot}}

\DeclareMathOperator{\Int}{Int}      
\DeclareMathOperator{\Cl}{Cl}        

\newcommand\iffdef{\;\mathrel{\mathord{:}\mathord{\longleftrightarrow}}\;}

\newcommand\iffslim{\longleftrightarrow}
\newcommand\defeq{\coloneqq}

\newcommand\rarrow{\longrightarrow}

\newcommand\fun{\ensuremath{\rightarrow}} 
\newcommand{\power}{\mathcal{P}}

\newcommand{\zero}{\mathsf{0}}
\newcommand{\one}{\mathsf{1}}

\newcommand{\Stone}{\mathtt{Ult}}

\DeclareMathOperator{\con}{\mathsf{C}} 

\newcommand\mathbackslash{\raisebox{.4pt}{\texttt{/}}}
\def\notcon{
  \renewcommand\stacktype{L}\mathrel{\ensurestackMath{%
  \ThisStyle{\stackon[0pt]{\SavedStyle\con}{\SavedStyle\mathbackslash}}}}%
}

\newcommand{\conf}{\mathbin{\dwakola}}

\def\nconf{%
\renewcommand\stacktype{L}\mathbin{\ensurestackMath{%
  \ThisStyle{\stackon[0pt]{\SavedStyle\dwakola}{\SavedStyle/}}}}}

\def\dwakola{
\mathord{\circ}\!\mathord{\circ}}

\newcommand{\prePt}{\mathbf{Q}_G}
\newcommand{\prePtW}{\mathbf{Q}_W}
\newcommand{\Abs}{\mathbf{A}}

\let\Grz=\Pt
\newcommand\Eq{\mathbf{Eq}}
\newcommand\Wthd{\mathbf{W}}

\newcommand{\topo}{\mathscr{O}}

\DeclareMathOperator{\RO}{RO} 

\newcommand{\fil}{\mathscr{F}}
\newcommand{\ult}{\mathscr{U}}

\newcommand{\End}{\mathbf{End}}

\newcommand{\Ult}{\mathbf{Ult}}
\newcommand{\Basis}{\mathscr{B}}

\newfont{\eurxi}{eurm10 at 10.95pt}
\newfont{\eurviii}{eurm7 at 8pt}
\newfont{\eurvii}{eurm7}

\newcommand{\Nat}{\mbox{\eurxi\symbol{33}}} 
\newcommand\Real{\mathds{R}} 

\def\dywiz{\kern0sp\discretionary{-}{}{-}\penalty10000\hskip0sp\relax}

\newcommand\BA{\boldsymbol{\mathsf{BA}}} 
\newcommand\GCA{\boldsymbol{\mathsf{GCA}}} 
\newcommand\Conc{\boldsymbol{\mathsf{Conc}}} 

\newcommand\Even{\mathds{E}}
\newcommand\Odd{\mathds{O}}

\DeclareSymbolFont{stx}{OMS}{txsy}{m}{n}
\DeclareMathSymbol{\Ing}{\mathrel}{stx}{118}
\DeclareMathSymbol{\ll}{\mathrel}{stx}{28}

\DeclareMathSymbol{\medcirc}{\mathord}{symbolsC}{7}
\DeclareMathSymbol{\Ext}{\mathrel}{symbolsC}{78}

\DeclareMathSymbol{\nIng}{\mathrel}{symbolsC}{64}
\DeclareMathSymbol{\Ov}{\mathrel}{symbolsC}{7}

\let\ingr=\Ing

\newcommand{\coloneqq}{\mathrel{\mathord{:}\mathord{=}}}

\newtheorem{theorem}{Theorem}

\newtheorem{corollary}{Corollary}

\theoremstyle{remark}

\newcommand\scrE{\mathscr{E}}
\newcommand{\scrG}{\mathscr{G}}

\renewcommand{\BA}{\boldsymbol{\mathsf{BA}}}
\newcommand{\DV}{\boldsymbol{\mathsf{DV}}}
\newcommand{\KHaus}{\boldsymbol{\mathsf{KHaus}}}
\renewcommand{\Stone}{\boldsymbol{\mathsf{Stone}}}
\DeclareMathOperator{\CO}{CO}
\renewcommand{\prePt}{\mathbf{Q}}
\newcommand\covers{\succeq}

\renewcommand{\Nat}{\omega}
\renewcommand{\End}{\mathbf{MRF}}

\title[Mathematical methods in region-based theories of space]{Mathematical methods in region-based theories of space:\\ the case of Whitehead points}

\author{Rafa\l\ Gruszczy\'nski}

\date{}
\address{Rafa\L\ Gruszczy\'nski\\
Department of Logic\\
Nicolaus Copernicus University in Toru\'n\\
Poland}

\email{gruszka@umk.pl}

\urladdr{www.umk.pl/\textasciitilde gruszka}

\begin{document}

\maketitle

\begin{abstract}
Regions-based theories of space aim---among others---to define \emph{points} in a geometrically appealing way. The most famous definition of this kind is probably due to \cite{Whitehead-PR}. However, to conclude that the objects defined are points indeed, one should show that they are points of a~geometrical or a topological space constructed in a~specific way. This paper intends to show how the development of mathematical tools allows showing that Whitehead's method of extensive abstraction provides a~construction of objects that are fundamental building blocks of specific topological spaces.

\medskip


\medskip

\end{abstract}

\section{Introduction}

Imagine yourself trying to read out space's structure from the flux of data that reach your senses. After \citet{Russell-OKEW} we might say that you are submerged in the perspective space---the space of private experience, a small fragment of the world. Yet your ambitions go way beyond that. You aim at a general mathematical theory that will reflect the essential, structural properties of a~large fragment of what we know as the universe. You know it is feasible. Faithful and efficient systems of geometry are exactly such theories, and they have been with us since antiquity. Developed at the outset as tools to handle practical problems of relatively small communities, they turned into theories describing universal properties of larger fragments of space, including the properties of the universe as such after the emergence of non-Euclidean geometries. The rise of topology has been driven by the search for space's most general features, as well as for the solution of real-world problems, Euler's \emph{K\"onigsberg bridges puzzle tour} to be one of them. Purely mathematical enterprise at the beginning, topology flourished as a~branch of mathematics with applications in macro- and micro-scale. All those achievements were obtained by experiencing fragments of our world only but turned out to be so powerful as to describe its most general properties.

Put yourself into the shoes of an admirer of geometry and topology who, at the same time, finds one thing to be a bit troubling---the fundamental constituents of geometrical and topological spaces are points, highly idealized, dimensionless objects that cannot be found in the space of private experience. Thus you ask yourself the question: \emph{could points be mathematically satisfactory explained employing the objects from the perspective space?}

One of the very first endeavors toward a~positive answer to the question was due to Alfred N. \citet{Whitehead-ECPHK,Whitehead-CN,Whitehead-PR}. He presented various constructions of points out of which the one from \emph{Process and reality} was best developed and gained the attention of the community of logicians, mathematicians, and philosophers.\footnote{To tell the truth, the constructions of \emph{points} from \citet{Whitehead-ECPHK,Whitehead-CN} were wrong, as observed by de Laguna \citep{deLaguna-EAAS}. The reason was that initially, Whitehead worked with \emph{part of} relation only, and de Laguna suggested---and rightly so---going beyond it and adopting the notion of \emph{containing} (the dual of the modern \emph{non-tangential inclusion}) as one of the primitives.} However, having defined points, the English mathematician never bothered himself to show that the entities constructed are building blocks of any space.

This paper's goals are very modest, as we aim to show how the development of formal methods from the XXth century let us carry out Whitehead's construction in a rigorous mathematical manner and formulate a~partial positive solution to the problem of existence of non-trivial topologies based on Whitehead points. In light of this, the paper does not provide any new groundbreaking results in the field of region-based topology but rather shows how various results obtained within it allow us to draw a positive conclusion concerning Whitehead points: not only do exist structures with Whitehead points, but these points are also building blocks of topological spaces that were constructed in the area of representation theory for Boolean algebras and their extensions.

\section{The informal construction}\label{sec:informal-construction}

Observe that the data you receive through your senses, concerning the spatial entities, contain various objects that we may collectively call \emph{regions}. Both the laptop on your desk and the courtyard you see from the window of your office are regions, chunks of space. Those chunks are related to each other in various ways, of which two seem to be the most general: (a) one region may be part of another, as the screen is part of the laptop, (b) two regions may touch each other, as in the case of the laptop and the surface of the desk, or can be separated, like the pen in your backpack and the cup of coffee standing next to your left hand. Next to these, we have the idea of points as precise locations in space. On the other hand, these can be represented as collections of shrinking regions in space, tapering down to the precise locations. One of the main driving forces of region-based theories is to capture this vague idea through parthood and contact.

\begin{figure}[!htbp]
    \centering
    \begin{center}
    \begin{tikzpicture}
    \foreach \y\z in {2.5/10,2.25/20,2/30,1.75/40,1.5/50,1.25/70,1/80,0.75/90,0.5/100}{
    \draw[fill=cyan!\z] plot[domain=0:350, smooth cycle] (\x:\y+rnd*0.2);
    }
    \foreach \x in {0,72,...,288}
    {
    \draw [dotted,shorten <= 2pt] (0,0) -- (\x:0.5cm);
    }
    \end{tikzpicture}
    \end{center}
    \caption{Point as a limit of shrinking system of regions}
    \label{fig:point-as-a-limit}
\end{figure}
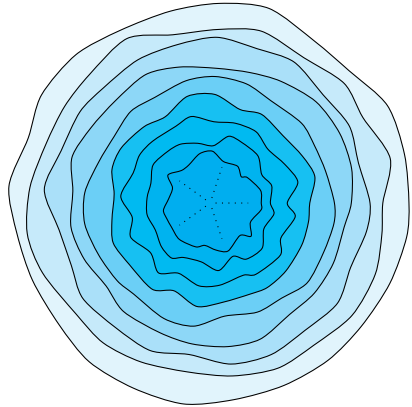

One way is to write down axioms that could be justified by how we seem to experience regions and their relations. We may engage both \emph{parthood} and \emph{contact}, as many authors did, or only just one of them, as was done originally by \citet{Whitehead-PR}. 
Let `$\ingr$' and `$\con$' be the two symbols that denote, respectively, \emph{parthood} and (binary) \emph{contact}. We read `$x\ingr y$' as \emph{$x$ is part of $y$} and `$x\con y$, as \emph{$x$ is contact with $y$} or $x$ \emph{touches} $y$. The most reasonable axioms for the former are probably those for one of the possible systems of mereology\footnote{See e.g. \citep{Pietruszczak-M-eng,Pietruszczak-FTP} and \citep{Varzi-M} for expositions of various mereological theories.} that is a~faithful representation of spatio-temporal properties between regions and their parts. For contact, the standard axioms to be assumed are the following: every region is in contact with itself: $x\con x$, the contact is symmetrical: if $x\con y$, then $y\con x$; if $x\con y$ and $y$ is part of $z$, then $x\con z$, which intuitively means that if $x$ touches $y$, then every region of which $y$ is part must also touch $x$. This is the axiomatic basis. Other axioms may be introduced, and we will get back to this in the sequel. The way of proceeding directly from the sense data to axioms of a~theory can be named, after \cite{Pratt-Hartmann-EARIRBTOS}, the \emph{empiricist} approach.

The objective is to capture the notion of \emph{point}, and the most famous such a try-out from the beginnings of region-based theories goes back to Alfred N. \citet{Whitehead-PR}. His complicated characterization of contact and definition of point extends over six pages of \emph{Process and reality} and is preceded by 24 assumptions and 15 other definitions, a solid overkill, to say the least. Let's get straight to the bottom of Whitehead's points as easy as it gets without delving into Whitehead's philosophical motivations. For these, we refer the interested reader to the excellent exposition of \cite{Varzi-PHOC}.

\begin{figure}[!htbp]
\begin{center}
\begin{tikzpicture}
\draw[fill=lightgray!50] plot[domain=0:350, smooth cycle] (\x:2+rnd*0.3) node at (1.7,0) {$x$};
\draw[fill=cyan] plot[domain=0:350, smooth cycle] (\x:1+rnd*0.3) node at (0,0)  {$y$};
\end{tikzpicture}
\end{center}
\caption{Region $y$ is a non-tangential part of region $x$.}\label{fig:ll}
\end{figure}
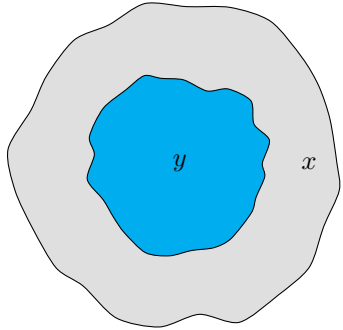

\citet{Whitehead-PR} follows the idea of the point from Figure~\ref{fig:point-as-a-limit}. To do this properly, one must first say what it means for one region to be \emph{a non-tangential part} (we will often use the phrase `well-inside' as a synonym of `non-tangential inclusion') of another: it is the case when the former is not in contact with the complement of the latter\footnote{If we are working in the classical mereology we have to be careful what we mean by the \emph{complement} as the zero region is absent. See \cite{Pietruszczak-M-eng} for details. In the case the main theory does not assume a~region that is the largest region, the notion of the complement may have no sense at all, and we have to define the situation from Figure~\ref{fig:ll} in a~different way. This can be done, e.g. by requiring that $y$ does not touch any region outside $x$. We refer the reader again to the paper by \cite{Varzi-PHOC}.} or, as we will often say, \emph{is separated from} the complement (see Figure~\ref{fig:ll}).

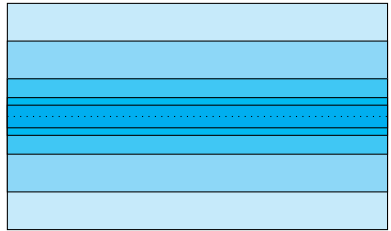
\begin{figure}[!htbp]
\begin{center}
\begin{tikzpicture}
\foreach \x\y\z in {0/3/20,0.5/2.5/40,1/2/60,1.25/1.75/80,1.35/1.65/100}
    {
        \draw[fill=cyan!\z] (0,\x) rectangle (5,\y);
    }
\draw[dotted] (0,1.5) -- (5,1.5);
\end{tikzpicture}
\end{center}\caption{A small fragment of an abstractive set of two-dimensional rectangles representing a one-dimensional segment, represented above by dots.}\label{fig:segment}
\end{figure}

On the Whiteheadian road to points, we begin with the definition of \emph{abstractive} sets of regions, i.e., these sets that:
\begin{enumerate}[label=(\alph*), align=left, leftmargin=\parindent,itemsep=0pt,leftmargin=*]
    \item do not any have minimum with respect to \emph{part of} relation, that is every region constituting an abstractive set has a proper part that is also in the set,
    \item their any two distinct elements are comparable with respect to \emph{non-tangential part} relation.
\end{enumerate}
The idea is that abstractive sets represent objects such as two-dimensional figures, planes, one-dimensional lines or segments (see Figure~\ref{fig:segment}), and---last but not least---points, as the readers will convince themselves looking at Figure~\ref{fig:point-as-a-limit} again. The question is \emph{how} to identify these abstractive sets that represent points? To this end, we define the \emph{covering} relation between abstractive sets as follows: $A$ \emph{covers} $B$ (in symbols: $A\succeq B$) iff for every region $x$ in $A$ there is a~region $y$ in $B$ such that $y$ is part of $x$. Now, if $A$ covers $B$ and vice versa, both sets represent the same object, and we can say those sets are equivalent. It is routine to verify that the equivalence of abstractive sets is indeed an equivalence relation: reflexive, symmetrical, and transitive. An equivalence class, say $[A]$, represents a~unique object and therefore deserves to be called \emph{a geometrical object}. Still, it does not have to satisfy our intuition of \emph{point} as dimensionless, <<infinitely>> small entity. How to identify these geometrical objects that do? A way out is via comparing geometrical objects in the following manner: $[A]\unrhd[B]$ if and only if $A\succeq B$. The relation $\unrhd$ is a~partial order, and if the partially ordered set of all geometrical elements happens to have minimal elements, then these elements truly deserve the name of \emph{points}.

Yet do they? How can we be sure that these are good candidates for points? After all, we have nothing to support this claim except for our intuition: when we \emph{think} about regions as extended objects of the spatio-temporal continuum, ordered by the aforementioned armory of relational concepts, then what we defined as points are abstract objects that are, in a way, so <<tiny>> that they must be good representations of what we may ever want to declare points. So far, so good, the problem is that the intuition may fail, and the best way to avoid failure is to put it to strict mathematical tests. To do this, we need proper \emph{formal} machinery, and thus we have to leave empiricism behind and take the path of \emph{rationalism}, as characterized by \cite{Pratt-Hartmann-EARIRBTOS}.

\section{The cornerstone}

How can we test objects for <<pointhood>>? The best method we have is that of the representation theory known from universal algebra, which allows us to show that given objects from some abstract or concrete algebraic structure are indeed points. The idea of \emph{representation} is a~formal embodiment of reducing the \emph{unfamiliar} and \emph{abstruse} to \emph{familiar} and \emph{comprehensible}. Or turning \emph{abstract} into \emph{concrete}.


To get the feeling of the mechanism of representation, let us divide the class of Boolean algebras into two subclasses of (a) \emph{abstract} and (b) \emph{concrete} Boolean algebras, respectively. Abstract BAs are defined as specific structures with a~distinguished domain whose elements are to satisfy certain conditions (axioms), usually for the binary operations of meet and join, the unary operation of complement, and the two individual constants: zero and one (unity). In the case of concrete BAs, we have a~fixed set $X$, and take as the domain of the algebra a~family $\mathcal S$ of its subsets that contains $X$ and $\emptyset$, and is closed for the set-theoretical operations of intersection, union, and complement. Such a~family is called a~\emph{field of sets}. More precisely, a~concrete BA may be identified with a~pair $\langle X,\mathcal S\rangle$ such that $\mathcal S$ is a~field of sets over~$X$ (see \citealp{Tao-TSALSRT}). It is evident that every concrete BA is an abstract one. It is also true that every abstract algebra is isomorphic to a~concrete algebra, although this statement is far from obvious. It was proved by Marshall \citet{Stone-TRBA}, who created the representation method relevant from this paper's point of view.

Stone's work's motivations were purely mathematical---he aimed to understand what Boolean algebras are and how they relate to other mathematical entities. The first step towards understanding was to show that given any~(abstract) Boolean algebra $B$ we can construct (in a~canonical way) a~concrete algebra $\langle X,\mathcal S\rangle$ that is isomorphic to~$B$.

With an algebra $B$ at hand, everything we have at our disposal is this algebra (plus various mathematical tools that are normally used). The situation is analogous to a~construction of a~term model of a~first-order theory by means of the Henkin method---we start with syntactical data, and we turn it into a~model of the theory. To tell the truth, Henkin's construction may be viewed as a~special case of the Stone theorem (see e.g., Exercises 4, 5, and 6 in pages 37--38 of \citealp{Koppelberg-HBA})

With every---either abstract or concrete---algebra, there is associated the notion of a~\emph{filter}, a~non-empty subset $\fil$ of (the domain of) $B$ that does not contain the zero element, is upward closed (in the sense of the standard Boolean order), and is closed for the binary meet operation.


The special place in the representation theory is occupied by \emph{ultrafilters}, i.e., filters that are maximal in the family of all filters (in the sense of set-theoretical inclusion), or equivalently, filters $\fil$ that satisfy the following condition: for any $x\in B$, either $x$ is in $\fil$ or its Boolean complement $-x$ is in~$\fil$. Given an algebra $B$, we will denote the family of all its ultrafilters by $\Ult(B)$, and we are going to use the letter `$\ult$' as ranging over $\Ult(B)$.

Applying set-theoretical machinery, we can prove that every non-zero object $x\in B$ is an element of an ultrafilter. Moreover, we can show that if $x$ and $y$ are distinct, then there is an ultrafilter $\ult$ that contains exactly one of these. Therefore, every object in $B$ can be unequivocally represented by all these ultrafilters to which it belongs. More formally, with every Boolean algebra $B$ we may associate an operation $\ult\colon B\to\power(\Ult B)$ (from the domain to the power set of the family of all ultrafilters of $B$) such that $\ult(x)\defeq\{\ult\in\Ult(B)\mid x\in\ult\}$ (\emph{the Stone mapping}) that is injective. It is routine to verify that the image of this operation:
\[
\ult[B]=\{\ult(x)\mid x\in B\}
\]
is a~field of sets. Indeed, $\Ult B=\ult(\one)\in\ult[B]$ and $\emptyset=\ult(\zero)\in\ult[B]$. $\ult[B]$ is closed for intersections and unions since:
\[
\ult(x)\cap\ult(y)=\ult(x\cdot y)\quad\text{and}\quad\ult(x)\cup\ult(y)=\ult(x+y)\,,
\]
where $\cdot$ and $+$ are the Boolean operations of meet and join, respectively; and the closure for set-theoretical complementation stems from the following equivalence:
\begin{equation}\label{eq:complement-for-Stone}
\ult\notin\ult(x)\iffslim x\notin\ult\iffslim -x\in\ult\iffslim\ult\in\ult(-x)\,.
\end{equation}

To conclude, to an abstract Boolean algebra $B$ we can always associate a~concrete isomorphic algebra $\langle\Ult(B),\mathcal S\rangle$ (with $\mathcal S\defeq\ult[B]$), that is isomorphic with~$B$, its canonical representation. This is the content of the set-theoretical version of the Stone representation theorem.

However, the construction may be carried on to a~topological representation. The main advantage of this is that it allows using spatial intuitions to draw consequences about the algebraic properties of Boolean algebras. In the case of set-theoretical representation above, the algebra $B$ is shown to be isomorphic to a~field of sets. In the case of the topological representation, it is proven that the field consists of distinguished---in one way or another---subsets of a~topological space.

With respect to these, two crucial observations are that (a) we may treat ultrafilters as points---building blocks of point-based topologies, (b) with the topological structure induced by sets $\ult(x)$ taken as basic open sets. The fact that $\ult[B]$ satisfies the conditions of a~basis stems from earlier observations for this family: every ultrafilter is in $\ult(\one)$, and $\ult(x)\cap\ult(y)=\ult(x\cdot y)$. Let $\mathscr S$ be the topology on $\Ult B$ with $\ult[B]$ as a~basis\footnote{Recall that a basis for a topology on the set $X$ is a family $\Basis$ of subsets of $X$ such that $X=\bigcup\Basis$ and for every $B_1,B_2\in\Basis$ and every $x\in B_1\cap B_2$ there is $B_3\in\Basis$ such that $x\in B_3\subseteq B_1\cap B_2$.}, the \emph{Stone topology}. The pair $\langle\Ult(B),\mathscr S\rangle$ bears the name of the \emph{Stone space} for the algebra~$B$.


Let us have a look at the basic features of Stone spaces. Firstly, observe that given any open basic set $\ult(x)$, it is a~straightforward consequence of \eqref{eq:complement-for-Stone} that its complement is open too. This means that the basis for $\mathscr S$ is built out of sets that are both closed and open (and are called \emph{clopen} for this reason). Such spaces are called \emph{zero-dimensional}, and they are not very intuitive from the point of view of properties of the perspective space. If we take, e.g., the three-dimensional Cartesian space that serves as the standard model of the (static) world around us, then we only find two clopen sets: the whole space $\Real^3$ and the empty set.
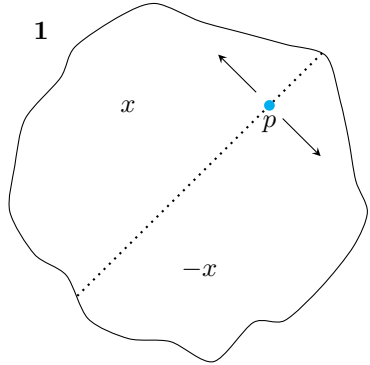
\begin{figure}[!htbp]
\begin{center}
\begin{tikzpicture}
    \draw[white,name path = line] (-2,-2) -- (2,2);
    \draw[name path = space] plot[domain=0:350, smooth cycle] (\x:2+rnd*0.5);
    \node at (-2,2) {$\mathbf{1}$};
    \draw[dotted,thick,name intersections={of=space and line}] (intersection-1) -- node [above left = 1cm] {$x$} node [below = 1cm] {$-x$} (intersection-2);
    \node[fill=cyan,circle,minimum size=4pt,inner sep=0pt,outer sep=5pt] (p) at (1,1) {} node at (p.south) {$p$};
    \draw[-stealth] (p.south east) -- ++(0.5,-.5);
    \draw[-stealth] (p.north west) -- ++(-0.5,.5);
\end{tikzpicture}
\end{center}\caption{In Stone spaces, points cannot be located on boundaries between regions, as there are no boundaries. The point $p$ is either a point of $x$ or a point of the complement of $x$.}\label{fig:border-point-1}
\end{figure}
For the other crucial property of Stone spaces, look at Figure~\ref{fig:border-point-1}. The intuition from the perspective space is that when we divide a~region into two parts, there is such a~thing as the boundary between the parts, and there are points that are located on the boundary. However, this is impossible in Stone spaces. The point $p$ from the figure is an~ultrafilter. Therefore either~$x$ is in $p$, or the Boolean complement of $x$ is in~$p$. In topological parlance, we say that the space is disconnected. For Stone spaces, the discontinuity phenomenon takes an extreme form: the only connected components of those spaces are singletons of points, i.e., the spaces are \emph{totally disconnected}. Again, this is not a~very intuitive property from the point of view of the perspective space. Actually, for the class of compact and Hausdorff spaces (a~larger class than that of Stone space), the two properties are equivalent, in the sense that every compact Hausdorff space is zero-dimensional iff it is totally disconnected (see Theorem 7.5 in \citealp[p.\,97]{Koppelberg-HBA}).

The aforementioned compactness is---in a way---a~topological version of finiteness: a space $X$ is \emph{compact} if for every family of open sets that covers the whole space, there is its finite subfamily that covers $X$ either. As every open set is the sum of some family of basic open sets, we may replace `open' with `basic open' in the definition. In the case of Stone spaces, compactness is a~consequence of the Ultrafilter Theorem, which says that every set $F$ of elements of a BA such that $F$ has the \emph{finite intersection property} is contained in an ultrafilter, where $F$ has the finite intersection property iff any finite subcollection of~$F$ has the non-zero meet: if $x_1,\ldots x_n\in F$, then $x_1\cdot\ldots\cdot x_n\neq\zero$. Using this, it is relatively easy to show that the Stone space $\Ult(B)$ is compact.


Another key feature of Stone spaces is the Hausdorff separation axiom: any two distinct points $x$ and $y$ can be separated by open sets, in the sense that there are disjoint open set $U$ and $V$ around $x$ and $y$, respectively. If ultrafilters $\ult_1$ and $\ult_2$ are distinct, there must be an~$x$ which is in only one of them, say $\ult_1$. But then $-x$ must be in $\ult_2$, and thus $\ult(x)$ and $\ult(-x)$ are
disjoint (basic) open sets around the two ultrafilters, i.e., points of the~Stone space.

To conclude, with every Boolean algebra $B$, we can associate a~topological space, the Stone space of $B$, which is Hausdorff, compact, and zero-dimensional.\footnote{Topological spaces that have these three properties are often called \emph{Boolean spaces}, and the name is used with the intention to treat such spaces somewhat independently from the Stone spaces of ultrafilters. However, as we will see, every Boolean space $X$ is a~Stone space, in the sense that we can associate with $X$ a~Boolean algebra $B$ whose Stone space $\Ult(B)$ is an exact copy of~$X$.} Moreover, the algebra $B$ is isomorphic to the family $\CO(\Ult(B)))$ of clopen sets of this space. Thusly, there is a~way from Boolean algebras to topological Stone spaces, i.e., structures with certain spatial data.

However, there is also a~way in the other direction. Any topological space $X$ carries a~Boolean algebra $\CO(X)$ of all its clopen subsets. In the case of Euclidean spaces $\Real^n$, this algebra will have only two elements: the whole space $\Real^n$ and the empty set. More generally, every connected space will carry the two-element Boolean algebra of its clopen subsets. Things get interesting if we limit our attention to Stone spaces only. In such a case, we obtain a~deep dependence between the class $\Stone$ of all Hausdorff, compact, and zero-dimensional spaces, and the class $\BA$ of all Boolean algebras.

Let us start with a Boolean algebra $B$. As we have seen, there is a~topological space that can be naturally associated with~$B$, the Stone space $\Ult(B)$. This space, on the other hand, carries a~Boolean algebra of its clopen subsets $\CO(\Ult(B))$, that is isomorphic to~$B$, i.e., $B$ and $\CO(\Ult(B))$ cannot be structurally distinguished.

\begin{figure}[!htbp]
\begin{center}
\begin{tikzcd}[row sep=large,column sep=large]
\BA\arrow[r, symbol=\ni] &[-3em]B \arrow[r, ->,dashed] \arrow[rd, "i" left=5pt, ->,shift right] \arrow[rd, "i^{-1}" right=5pt, <-,shift left] & \Ult(B) \arrow[d, dashed] \arrow[r, symbol=\in]& [-3em] \Stone\\
&{}& {\CO(\Ult(B))}\arrow[r, symbol=\in]& [-3em] \BA
\end{tikzcd}\caption{Any Boolean algebra $B$ is indistinguishable from the Boolean of clopen sets of the Stone space of $B$.}
\end{center}
\end{figure}
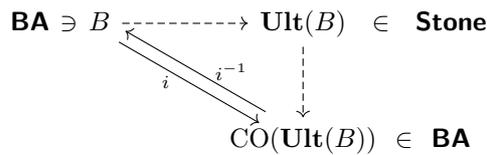

On the other hand, if we start from a~Stone space $X$, then $\CO(X)$ is a~Boolean algebra, and $\Ult(\CO(X))$ is its Stone space, that is, as the reader could expect, indistinguishable (\emph{homeomorphic} is the technical jargon) from~$X$.

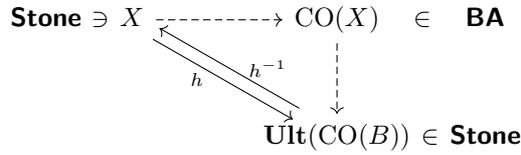
\begin{figure}[!htbp]
\begin{center}
\begin{tikzcd}[row sep=large,column sep=large]
\Stone\arrow[r, symbol=\ni] &[-3em]X \arrow[r, ->,dashed] \arrow[rd, "h" left=5pt, ->,shift right] \arrow[rd, "h^{-1}" right=5pt, <-,shift left] & \CO(X) \arrow[d, dashed] \arrow[r, symbol=\in]& [-3em] \BA\\
&{}& {\Ult(\CO(B))}\arrow[r, symbol=\in]& [-3em] \Stone
\end{tikzcd}\caption{Any Stone space $B$ is indistinguishable from the Stone space of the Boolean algebra $\CO(X)$ of the clopen sets of $X$.}
\end{center}
\end{figure}

Thus, \citet{Stone-TRBA} demonstrated that there is a~kinship between the world of algebraic structures and the world of topological spaces. In particular, when we focus on Boolean algebras and Stones spaces, the bondage is so strong that we can say they are two sides of the same coin or two aspects of the same abstract phenomenon.\footnote{The kinship also extends to homomorphisms between algebras and continuous mappings between the spaces, in the sense that to every homomorphism between BAs corresponds a~continuous mapping between their Stone mapping, and vice versa---with every continuous mapping between Stone spaces, there is associated a homomorphism between the algebras of their clopen sets. It is, roughly, the content of the famous Stone duality between the categories of Boolean algebras with homomorphisms, and Stone spaces with continuous mappings. For details, see, e.g., \citep{Johnstone-SS}.}

\begin{figure}[!htbp]
\begin{center}
\begin{tikzpicture}[node distance = 3cm,on grid]

\node[circle,draw,outer sep=2pt, minimum size = 2cm] (algebra) {algebra};
\node[circle, draw, outer sep=2pt, minimum size = 2cm, right = of algebra] (topology) {topology};

\path [-stealth, thick]
	(algebra.north east) edge [bend left] node [above] (F) {F} (topology.north west)
	(topology.south west) edge [bend left] node [below] (G) {G} (algebra.south east);

\node [below left = of algebra] (Boole) {Boolean algebras};
\node [above right = of topology] (Stone) {Stone spaces};

\node [above left = of F] (ult) {$\Ult$};
\node [below right = of G] (co) {$\CO$};

\path [-stealth, thick,cyan] (Stone.south) edge [bend left] (topology.east);
\path [-stealth, thick,yellow] (Boole.north) edge [bend left] (algebra.west);
\path [-stealth,thick,yellow] (ult.east) edge [bend left] (F.north);
\path [-stealth,thick,cyan] (co.west) edge [bend left] (G.south);


\end{tikzpicture}\caption{Boolean algebras and Stone topological spaces are very closely related.}\label{fig:alg-top}
\end{center}
\end{figure}
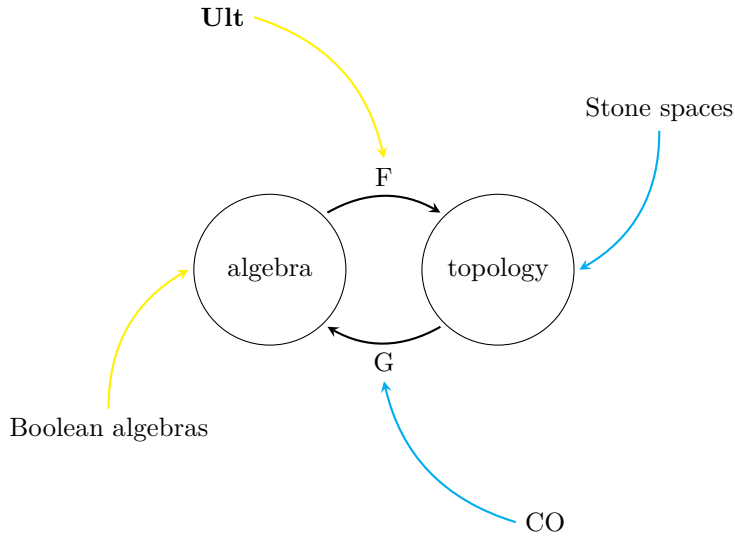

\section{The extension: De Vries algebras}

If we look at the diagram in Figure~\ref{fig:alg-top} we see that there are two very natural ways towards extending Stone results: we may encompass a larger (or a~different) class of topological spaces, but we may also tinker with algebras taken into account. One of such tinkerings may, in particular, involve extending the signature (the non-logical language).

\newlength{\rozszerzenie}
\settowidth{\rozszerzenie}{rozszerzenie}

\begin{figure}[!htbp]
\begin{center}
    \begin{tikzpicture}
    \tikzset{
                every node/.style={draw,rounded corners}
            }
        \node (twoways) {Two ways};
        \node [below left = of twoways] (general) {\parbox{\rozszerzenie}{\centering\small extension of the class of spaces}};
        \node [below right = of twoways] (specific) {\parbox{\rozszerzenie}{\centering\small extension of the signature}};
        \node [below = of general] (frames) {\parbox{\rozszerzenie}{\centering\small frames and locales}};
        \node [below right = of specific] (bcas) {\parbox{\rozszerzenie}{\centering\small contact algebras}};
         \node [below left = of specific] (devries) {\parbox{\rozszerzenie}{\centering\small De Vries algebras}};
        \draw[myarrow] (twoways.south west) -- (general.north east);
        \draw[myarrow] (twoways.south east) -- (specific.north west);
        \draw[myarrow] (general.south) -- (frames.north);
        \draw[myarrow] (specific.south east) -- (bcas.north west);
        \draw[myarrow] (specific.south west) -- (devries.north east);
    \end{tikzpicture}\caption{Possible extensions of the Stone duality}
    \end{center}
\end{figure}
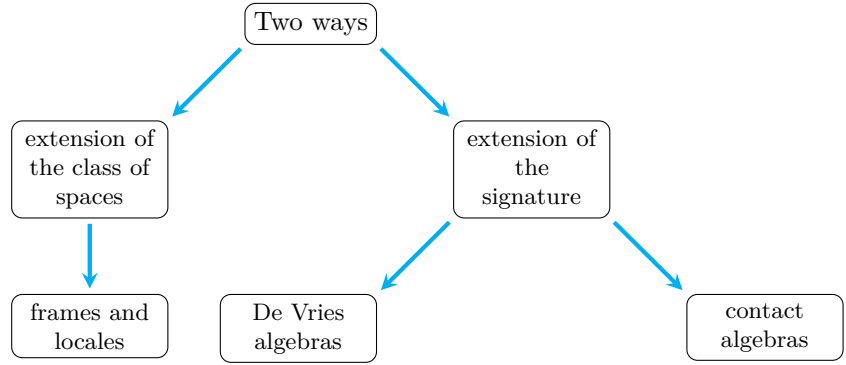

The extension of the class of topological spaces leads to a~fruitful and fascinating theory of frames and locales (see, e.g., \citealp{Johnstone-SS}, \citealp{Picado-Pultr-FL,Picado-Pultr-SPFT}), which in a~nutshell can be described as a~region-based theory of space in which the notion of \emph{open set} is taken as basic. It is probably most developed among all region-based approaches. Yet, its objectives and main motivations (for the exposition of these see, e.g., \citealp{Johnstone-PPT}) are not, at least directly, connected to the leading topic of this paper.\footnote{\cite{Mormann-CLWTS} presents a~solution of what he calls \emph{a Whitehead's problem} in the framework of Heyting algebras and continuous lattices, structures that are of particular importance in the theory of frames and locales. However, his paper does not mention Whitehead points and instead constructs topological spaces whose points are \emph{Dedekind ideals} (in the terminology adopted by us further in the paper, these could be called \emph{round} ideals). This is because \citeauthor{Mormann-CLWTS} defines the Whitehead problem as constructing spaces of points from regions of a~uniform dimension that sets of points can faithfully represent. If the reader wishes, they may think about our paper as presenting a~solution to the same problem yet utilizing the specific technique of conjuring up points put forward in \emph{Process and reality}.}  This is mainly due to the fact that we want to narrow down the notion of \emph{region} to those interpretations that are faithful models of fragments of the perspective space, while the notion of \emph{open set} is probably the most encompassing among primitive concepts of point-free theories.

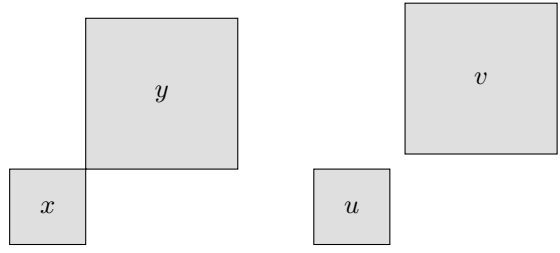
\begin{figure}[!htbp]
    \begin{center}
        \begin{tikzpicture}
\draw[fill=lightgray!50] (0,0) rectangle node {$x$} (1,1) rectangle node {$y$} (3,3);
\draw[xshift=4cm,fill=lightgray!50] (0,0) rectangle node {$u$} (1,1) (1.2,1.2) rectangle node {$v$} (3.2,3.2);
\end{tikzpicture}
        \caption{The regions $x$ and $y$ are incompatible and touch one other, while $u$ and $v$ are incompatible and separated.}
        \label{fig:touch-separate}
        \end{center}
    \end{figure}

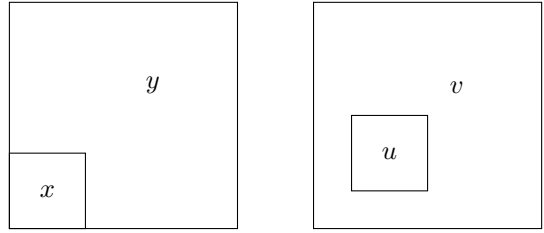
\begin{figure}[!htbp]
    \begin{center}
        \begin{tikzpicture}
\draw (0,0) rectangle node {$x$} (1,1);
\draw (0,0) rectangle node[above right,outer sep = 5pt] {$y$} (3,3);
\draw[xshift=4cm] (0.5,0.5) rectangle node {$u$} (1.5,1.5);
\draw[xshift=4cm] (0,0) rectangle node[above right,outer sep = 5pt] {$v$} (3,3);
\end{tikzpicture}
        \caption{The region $x$ is part of $y$ but touches the complement of $y$, while $u$ is well-inside $v$.}
        \label{fig:below-waybelow}
        \end{center}
    \end{figure}

Thus, we will follow the way of the signature extensions. The reasons to do this may vary, and our primary motivation is that the language of Boolean algebras does not differentiate between situations in which regions are incompatible (in the sense that their Boolean product is zero) and separated, and those where regions are incompatible but touch each other (see Figure~\ref{fig:touch-separate}). Equivalently, Boolean algebras cannot discern the difference between the situation in which $x$ is part of $y$ but does not touch the complement of $y$, and the one in which $x$ is part of $y$ and touches the complement of $y$, i.e., from the point of view of Boolean algebras there is no difference between the two scenarios in Figure~\ref{fig:below-waybelow}.

Again, there may be different reasons to ponder Boolean algebras' extensions, either with the touching relation or well-inside relation. As we already saw in Section~\ref{sec:informal-construction}, proper (whatever it means now) construction of points may require it. Nevertheless, the reasons may be less philosophical and more practical as in the case of \cite{deVries-CSC} work, which will serve as our starting point towards the justification of Whitehead's construction.

De Vries's aim was mathematical at heart: algebraization of the notion of \emph{compactness} of a~topological space. De Vries's algebras are just (complete) Boolean algebras extended with a~binary relation $\ll$ whose intended interpretation is non-tangential inclusion (well-inside part). The axioms concerning $\ll$ are the following:
\begin{gather*}
        \one\ll\one\tag{DV1}\,,\label{DV1}\\
        x\ll y\rarrow x\leq y\tag{DV2}\,,\label{DV2}\\
        x\leq y\wedge y\ll z\wedge z\leq w\rarrow x\ll w\tag{DV3}\,,\label{DV3}\\
        x\ll y\wedge x\ll z\rarrow x\ll y\cdot z\tag{DV4}\,,\label{DV4}\\
        x\ll y\rarrow -y\ll -x\tag{DV5}\,,\label{DV5}\\
        x\ll y\rarrow \exists z\,(x\ll z\wedge z\ll y)\tag{DV6}\,,\label{DV6}\\
        (\forall x\neq\zero)(\exists y\neq\zero)\,y\ll x\,.\tag{DV7}\label{DV7}\\
    \end{gather*}
These may not be self-evident at first sight, so let us explain them in a~proper setting. The concrete De Vries algebras can be obtained from regular open algebras of $\kappa$-normal topological spaces.\footnote{A space $X$ is \emph{$\kappa$-normal} (or \emph{weakly normal}) iff any pair of its disjoint regular closed sets can be separated by open sets (see \citealp{Shchepin-RFNNS}).} A subset $x$ of a~topological space is \emph{regular open} if $x$ is equal to the interior of its closure: $x=\Int\Cl x$.\footnote{Alternatively, regular open sets can be characterized as regular elements in the lattice $\Omega(X)$ of all open sets of $X$. Such a lattice is a~Heyting algebra and thus may have elements that are not regular, in the sense that if $x^\ast$ is a~relative complement of $x$, then $x^{\ast\ast}\nleq x$ (the reverse inclusion is always true). Thus $x$ is regular open if $x=x^{\ast\ast}$.} From a geometrical point of view, regular open sets of $\Real^n$ are those open sets that do not have <<surprises>> in the form of cracks, holes, punctures, or snags. For this reason, they are sometimes considered good candidates for mathematical regions of the perspective space.\footnote{Nowadays the class of all regular open sets of $\Real^n$ is usually considered too large to model regions of the surrounding world. Various authors put forward different limitations on it, see, e.g., \cite{DelPiero-CFRUSCM,DelPiero-NCFR}, \cite{Lando-et-al-ACORRBMAT}, \cite{Pratt-Lemon-OPPM}, and \cite{Schoop-PiPFM}.\label{foot:limitations}}

As is well known, given a~topological space $X$, the family of $\RO(X)$ of its regular opens is a~complete Boolean algebra with the operations defined as follows:
\begin{gather*}
U\cdot V\defeq U\cap V\qquad U+V\defeq\Int\Cl(U\cup V)\qquad -U\defeq\Int(X\setminus U)\\
\bigvee_{i\in I}U_i\defeq\Int\Cl\bigcup_{i\in I}U_i\qquad \bigwedge_{i\in I}U_i\defeq\Int\bigcap_{i\in I}U_i\,.
\end{gather*}

If we interpret the non-tangential inclusion in the standard way as:
\[
U\ll V\iffdef\Cl U\subseteq V\,,
\]
and assume that space $X$ is $\kappa$-normal, then we will see that $\RO(X)$ is a~De Vries algebra, with \eqref{DV6} (the so-called \emph{interpolation axiom}) corresponding to the $\kappa$-normality of the space, and \eqref{DV7} to its weak version of regularity according to which every non-empty open set $V$ has a~nonempty set $U$ whose closure is a~subset of~$V$. Since there are $\kappa$-normal spaces, there are De Vries algebras.

\begin{figure}[!htbp]
\begin{center}
\begin{tikzpicture}
\draw[name path = space] plot[domain=0:350, smooth cycle] (\x:2+rnd*0.3) node [, right, inner sep = 0.5cm] {$y$};
\draw[cyan] plot[domain=0:350, smooth cycle] (\x:1+rnd*0.3) node at (0,0)  {$x$};
\draw[red] plot[domain=0:350, smooth cycle] (\x:1.3+rnd*0.2) node [right]  {$z$};
\end{tikzpicture}\caption{A geometrical interpretation of \eqref{DV6} axiom: between any two regions $x$ and $y$ such that $x$ is well-inside $y$ we can squeeze in a~third region well-above $x$ and well-inside $y$.}
\end{center}
\end{figure}
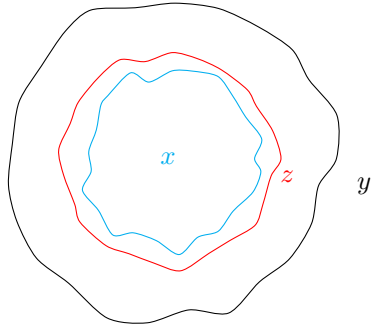

To pin down points in a De Vries algebra, the Stone-like technique of treating points as sets of regions is applied. To this end, the~family of \emph{round} filters is distinguished, i.e., filters $\fil$ that have the following property:
\[
(\forall x\in\fil)(\exists y\in\fil)\,y\ll x\,.
\]
It is easy to see that every De Vries algebra must have a~round filter: $\{\one\}$, trivial as it is. A less trivial example may be obtained if there is a non-zero region distinct from the unity, say~$x$. Then, by \eqref{DV7} and the Axiom of Dependent Choices, we can come up with a sequence of non-zero elements:
\[
\ldots x_2\ll x_1\ll x_0=x\,,
\]
and the filter generated by the sequence: $\fil\defeq\{y\mid(\exists{n\in\Nat})\,x_n\ll y\}$ must be round. An easy application of the Kuratowski-Zorn lemma shows that there exist \emph{maximal} round filters, and they are meant to be points of spaces.\footnote{The original \citet{deVries-CSC} terminology was different: he called \emph{concordant} and \emph{maximal concordant filters} round filters and maximal round filters, respectively. With other authors, the reader may also encounter terms \emph{contracting} and \emph{maximal contracting} filters. The latter are often called \emph{ends} in the framework of proximity approach the mereotopology. We have decided to use `round' as it is currently the most established practice among researchers within the field.}

Let us have a look at two concrete examples. Take the real line $\Real$, which is a~normal (and the more so $\kappa$-normal) space. Consider the family of intervals $\{(-1-\nicefrac{1}{n},1+\nicefrac{1}{n})\mid n\in\Nat^+\}$ whose elements are regular open in $\Real$. The filter $\fil$ generated by this family is round. However, it is not maximal in the family of round filters. We can extend $\{(-1-\nicefrac{1}{n},1+\nicefrac{1}{n})\mid n\in\Nat^+\}$ with some regular open sets well-inside $(-1,1)$ which will result in a~proper extension of $\fil$. Thus $\fil$ does not represent a point, which is good, since what all the regions in $\fil$ have in common is the interval $(-1,1)$, a~continuum of points. For a~positive example, take any point $x\in\Real$ and let $\RO(\topo_x)$ be the family of all regular opens around $x$. $\RO(\topo_x)$ is obviously a~filter, and since the real line is regular, it is round. But it must also be maximal. For suppose $\fil$ is a~round filter extending $\RO(\topo_x)$, and let $V$ be an element of~$\fil$ but not of $\RO(\topo_x)$. Therefore $x\notin V$. $\fil$ is round, so there is $M\in\fil$ with $\Cl M\subseteq V$. Thus, regularity entails existence of a regular open set $R$ around $x$ that is disjoint from $M$. But both $R$ and $M$ are elements of $\fil$, so $\fil$ is an improper filter (i.e. $\fil=\RO(X)$). Thus $\RO(\topo_x)$ is maximal round filter. This is good since the family of all regular open sets around~$x$ should uniquely determine~$x$.

Due to the latter example, we may be tempted to think that De Vries might have a  geometrical intuition of point similar to Whitehead's. However, a~certain example shows that the ideas of compactness and compactification were the leading ones for the Dutch mathematician, and it's a~point of discrepancy between his and Whitehead's approach. Consider\label{page:example-DV-not-W} the following chain of regions of $\RO(\Real)$: $\{(n,+\infty)\mid n\in\omega\}$. The filter $\fil$ that it generates is round and thus is contained in a~maximally round filter $\fil'$, a~point. This filter represents a~point at infinity in $\Real$, since it cannot be $\RO(\topo_x)$ for any real number~$x$. See also Figure~\ref{fig:point-at-infinity} for a~geometrical intuition in the case of two-dimensional space.

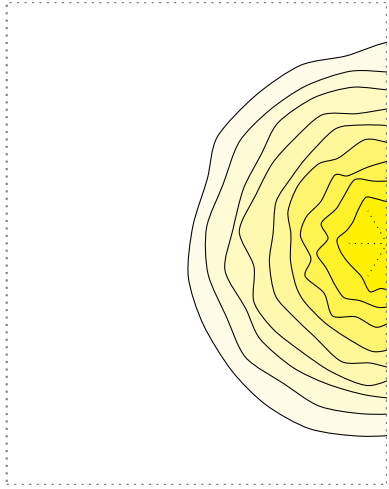
\begin{figure}[!htbp]
\begin{center}
    \begin{tikzpicture}
    \clip (0,-3.2) rectangle (-5.2,3.2);
    \foreach \y\z in {2.5/10,2.25/20,2/30,1.75/40,1.5/50,1.25/70,1/80,0.75/90,0.5/100}{
    \draw[name path = space,fill=yellow!\z] plot[domain=0:350, smooth cycle] (\x:\y+rnd*0.2);
    }
    \foreach \x in {120,180,240}
    {
    \draw [dotted,shorten <= 2pt] (0,0) -- (\x:0.5cm);
    }
    \draw[dotted,thick,gray] (0,-3.2) rectangle (-5,3.2);
    \end{tikzpicture}\caption{De Vries points involve points at infinity.}\label{fig:point-at-infinity}
    \end{center}
\end{figure}

Why do we maintain that this example shows that Whitehead and De Vries had different objectives?  The thing is that if we are to treat points as unique locations in the perspective space, points at infinity do not fit into this. Figuratively speaking, they are too far from our experience to enter the domain of points. At the very end of Section~\ref{sec:solution} we will demonstrate that the above constructed maximal contracting filter is not a~point in the sense of Whitehead.

Maximal round filters are exactly those round filters that satisfy the following condition:
\begin{equation}\tag{$\dagger$}\label{eq:ends-conditions}
    x\ll y\rarrow-x\in\fil\vee y\in\fil\,,
\end{equation}
or, as the readers may easily convince themselves:
\begin{equation}\tag{$\ddagger$}\label{eq:ends-conditions-2}
    x\notcon y\rarrow -x\in\fil\vee -y\in\fil\,,
\end{equation}
where contact is defined as:
\[
x\con y\iffdef \neg(x\ll -y)\,.
\]
On the other hand, for the whole class of filters of a Boolean algebra we have that $\fil$ (not necessarily round) is an ultrafilter if and only if:
\begin{equation*}
    x\cdot y=\zero\rarrow -x\in\fil\vee -y\in\fil\,.
\end{equation*}
Therefore, if we have additional information that $\ll$ coincides with $\leq$ (or, equivalently, contact is overlap), the family of maximal round filters is exactly the family of ultrafilters since every region is incompatible with its complement. However, in general, we cannot exclude existence of points living on the borders of regions and their complements, as we did in the case of spaces of ultrafilters (see figures~\ref{fig:border-point-1} and~\ref{fig:border-point-2}). Even more can be said: if $x$ is in contact with $-x$, then there is a~maximally round filter $\scrE$ such that $x\notin\scrE$ and $-x\notin\scrE$. This leads to an interesting conclusion: if every non-zero region is in contact with its complement, then the space of maximal round filters should be connected (if only there are such spaces).

\begin{figure}[!htbp]
\begin{center}
\begin{tikzpicture}
    \draw[white,name path = line] (-2,-2) -- (2,2);
    \draw[name path = space] plot[domain=0:350, smooth cycle] (\x:2+rnd*0.5);
    \node at (-2,2) {$\mathbf{1}$};
    \draw[dotted,thick,name intersections={of=space and line}] (intersection-1) -- node [above left = 1cm] {$x$} node [below = 1cm] {$-x$} (intersection-2);
    \fill[cyan] (1,1) circle (2pt) node [below,black] {$p$};
\end{tikzpicture}
\caption{In spaces of maximally round filters, points may inhabit the boundaries of regions and their complements, if the regions and their complements are in contact, in the sense that $\neg (x\ll x)$.}\label{fig:border-point-2}
\end{center}
\end{figure}
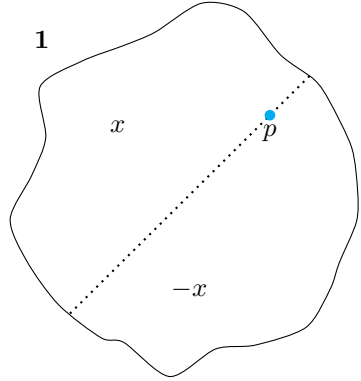

There are, of course. The standard Stone-like assignment $\scrE\colon B\fun\power(\End(B))$, where $B$ is a De Vries algebra and $\End(B)$ is the set of all its maximal round filters leads to the family $\Basis\defeq\{\scrE(x)\mid x\in B\}$ which satisfies the standard properties of a~basis. The spaces $\langle\End(B),\topo\rangle$ thus constructed are Hausdorff, since if $\scrE_1\neq\scrE_2$, and there is a~region $x$ in, say $\scrE_1\setminus\scrE_2$, then there is a~region $y\in\scrE_1$ well-inside~$x$. By \eqref{eq:ends-conditions} either $-y\in\scrE_2$ or $x\in\scrE_2$, and since the second disjunct does not hold, the first is true. But $\scrE(y)\cap\scrE(-y)=\emptyset$, and $\scrE_1\in\scrE(y)$ and $\scrE_2\in\scrE(-y)$.

Every $\End(B)$ must also be compact. The proof is slightly more complicated than the one for Stone spaces, and we skip it not to mar the gist of the paper with unnecessary technicalities. The important thing is that to every De Vries algebra corresponds a~certain topological space that is Hausdorff compact. Similarly to the situation for Boolean algebras and Stone spaces, given a~topological space $X$ that is Hausdorff compact, its family of $\RO(X)$ with $\ll$ interpreted as the topological well-inside inclusion must be a~De Vries algebra. Again, if we start with $B$, go to $\End(B)$ and to $\RO(\End(B))$, then we have that either $B$ can be densely embedded in $\RO(\End(B))$, or is isomorphic with $\RO(\End(B))$, if complete. Since we changed the class of algebras from Boolean to De Vries we need an appropriate notion of isomorphism that remains essentially the same as the standard one, with an extra condition stipulating that $\ll$ is preserved in the following sense:
\[
x\ll y\iffslim h(x)\ll h(y)\,.
\]
We will call the mapping $h$ \emph{De Vries isomorphism}. To repeat, if the initial algebra $B$ is not complete, then the mapping $\scrE\colon B\to\RO(\End(B))$ is a~dense De Vries embedding\footnote{The embedding $\scrE$ is dense in quite a~strong sense, that is if $x,y\in \RO(\End(B))$ are such that $x\ll y$, then there is $z\in B$ for which $x\ll \scrE(z)\ll y$.} of $B$, and in the case $B$ is complete, the same mapping is a De Vries isomorphism.

\begin{figure}[!htbp]
\begin{center}
\begin{tikzcd}[row sep=large,column sep=large]
\DV\arrow[r, symbol=\ni] &[-3em]B \arrow[r, ->,dashed] \arrow[rd, "\scrE" left=5pt, ->,shift right,hook] & \End(B) \arrow[d, dashed] \arrow[r, symbol=\in]& [-3em] \KHaus\\
&{}& {\RO(\End(B))}\arrow[r, symbol=\in]& [-3em] \DV^c
\end{tikzcd}\caption{Let $\DV$ be the class of De Vries algebras, and $\KHaus$ the class of compact topological spaces. Any De Vries algebra embeds densely into the algebra of regular open sets of the compact Hausdorff space for $B$.}\label{fig:DV-to-KHaus}
\end{center}
\end{figure}

\begin{figure}[!htbp]
\begin{center}
\begin{tikzcd}[row sep=large,column sep=large]
\DV^c\arrow[r, symbol=\ni] &[-3em]B \arrow[r, ->,dashed] \arrow[rd, "\scrE" left=5pt, ->,shift right] \arrow[rd, "\scrE^{-1}" right=5pt, <-,shift left] & \End(B) \arrow[d, dashed] \arrow[r, symbol=\in]& [-3em] \KHaus\\
&{}& {\RO(\End(B))}\arrow[r, symbol=\in]& [-3em] \DV^c
\end{tikzcd}\caption{Let $\DV^c$ be the class of complete De Vries algebras. Any its element $B$ is indistinguishable from the De Vries algebra of regular open sets of the compact Hausdorff space for $B$.}\label{fig:complete-DV-to-KHaus}
\end{center}
\end{figure}

\begin{figure}[!htbp]
\begin{center}
\begin{tikzcd}[row sep=large,column sep=large]
\KHaus\arrow[r, symbol=\ni] &[-3em]X \arrow[r, ->,dashed] \arrow[rd, "h" left=5pt, ->,shift right] \arrow[rd, "h^{-1}" right=5pt, <-,shift left] & \RO(X) \arrow[d, dashed] \arrow[r, symbol=\in]& [-3em] \DV^c\\
&{}& {\End(\RO(X))}\arrow[r, symbol=\in]& [-3em] \KHaus
\end{tikzcd}\caption{Any compact Hausdorff topological space $X$ is homeomorphic to the compact Hausdorff space of the the complete De Vries algebra of regular open sets of $X$.}\label{fig:KHaus-to-DV}
\end{center}
\end{figure}

To conclude, De Vries, through pursuing his algebraic objectives, showed a~way to represent structures with a~version of a~point-free topological nearness as fully-fledged topological spaces. In the next section, we will see how it helps to understand another classical point-free topology by a~Polish logician Andrzej Grzegorczyk, which on the other hand, will let us show that Whitehead points (or at least some of them) are indeed points of a~certain class of topological spaces.

\section{The criterion of points}

From the above, we can see that we have a general method of constructing spaces from algebraic data \emph{via} mimicking Stone's technique to treat points as subsets of the domain. So what would it mean to achieve the Whitehead's goal?, i.e., explain points in a~geometrically appealing way.  On the intuitive level, Whitehead's points are collections of regions related to each other via spatially motivated relations. The intuition may be turned into a precise notion in two steps: firstly,
by imposing an algebraic structure on regions to reflect the most general properties of the perspective space (i.e., extend the signature); secondly, by showing  that the Whitehead's minimal geometrical objects  reconstructed within such a structure as higher-order objects are indeed points of a~certain space.

More precisely---\label{page:method}and more generally---suppose $\langle A,R_1,\ldots,R_n\rangle$ is an algebraic structure with relations $R_1,\ldots,R_n$, all these together modelling the universe of regions. Suppose $\mathscr{P}$ is the set of higher-order objects defined within this structure. The main problem is now to find a~topological space with $\mathscr{P}$ as the underlying set of points (similarly as ultrafilters are taken as points of Stone spaces, and maximally round filters as points of compact Hausdorff spaces) that naturally models regions (elements of the domain) and relations $R_i$. That is, $\mathfrak{A}\defeq\langle A,R_1,\ldots,R_n\rangle$ is captured within $\mathscr{P}$ as $\mathfrak{U}'=\langle A',R_1',\ldots,R_n'\rangle$ in a similar way as any Boolean algebra $B$ is captured as the algebra of clopen sets of its Stone space $\Ult(B)$. This, in particular, means that $A'$ is a family of subsets of $\mathscr{P}$, and that the set $\mathscr{P}$ of points may be given an appropriate topology in which every $R_i$ can be modeled in such a~way that $R_i$ holds among regions iff $R'_i$ holds among their point-based counterparts. Usually, subsets of regular open (or regular closed) sets of $\mathscr{P}$ are taken as models of regions (see footnote~\ref{foot:limitations}). If this has been achieved, then we may say we have the solution to the problem of points, as we have the representation of the original structure in the structure built from higher-order elements of $\mathscr{P}$ that now deserve the name of \emph{points}. In this manner, every BA is represented in the space of ultrafilters, and every De Vries algebra in the space of maximally round filters. This justifies naming both ultrafilters and maximally round filters as points. The idea is now to repeat the above steps with a proper algebraic structure in place of $\mathfrak{A}$, and a~set of Whitehead points in place of $\mathscr{P}$.

In light of the theorems of Stone's and De Vries's, one could naturally ask could either ultrafilters of maximally round filters serve as Whitehead points? Why do they fall short? In the case of ultrafilters, the main problem is hidden in the fact that if they are points, the contact relation collapses to overlap. Indeed, suppose we have a mapping $f$ that represents regions of a~Boolean algebra in $\power(\Ult B)$, and that $f$ is the standard Stone-like function, that is, for every region $x$ its points are all these ultrafilters $\ult$ that has $x$ as an element. But then, as we observed earlier (see page \pageref{fig:border-point-1}), there are no points on the boundaries of the regions, so the contact can only be the overlap, i.e., we cannot model the situations in which objects are external to each other and touch each other at the same time. Moreover, if contact and overlap coincide, in complete algebras, there are no Whitehead points \citep{Gruszczynski-et-al-GAWPTSC}, so in general, ultrafilters cannot serve as them.

This does not mean that ultrafilters are always bad candidates for building blocks of spaces of points in which the contact relation is to be modeled. A sophisticated and elegant theory of this kind was created by Peter \cite{Roeper-RBT}. Yet his points are not ultrafilters themselves, but equivalence classes of ultrafilters that, in the end, can be shown to be maximal round filters \citep{Gruszczynski-NTP}.

As for De Vries points, we have shown above that their class is too large for the class of Whitehead points, in the sense that Whitehead points may only be a proper subset of the set of all maximal round filters. We'll get back to this problem in Section~\ref{sec:solution}.

For the completeness of the presentation, it should be emphasized that higher-order constructions are not the only method of explaining points, and some scholars either defined points in terms of regions (elements of the domain) or distinguished a~subset of the domain as the collection of points. The most important examples of this kind of approach---in either topological or geometrical setting---are \cite{Eschenbach-AMDOP}, \cite{Huntington-SPAGETSRI}, \cite{Galton-MM}, \cite{Galton-TMODS}, \cite{Hahmann-CODI}, \cite{Robering-TWIGTTP}, and \cite{Schoop-PiPFM}.

\section{Grzegorczyk points}\label{sec:Grzegorczyk}

This time we start with the contact relation. The difference is irrelevant from a~logical point of view, as with enough axioms the two approaches, either via contact or via non-tangential inclusion, are definitionally equivalent. However, the proper terminology well-chosen at the outset will equip us with a~user-friendly language. The purpose is to expose the point-free topology of Grzegorczyk's (\citeyear{Grzegorczyk-AGWP}), who by the way chose the third way and based his system on the notion of \emph{separation}, yet this is again an~equivalent approach to those used in this paper. \footnote{Strictly speaking, Grzegorczyk did not work with Boolean algebras, but with mereology, which is closely related to the former, see, e.g., \cite{Pietruszczak-M-eng}. The differences are mainly hidden in technical intricacies, as mereologies generally do not have zero elements and are thus semi-lattices.}

Before we begin a~proper, mathematical exposition of Grzegorczyk's construction and before we draw an analogy between this and Whitehead's, let us remind that it was \cite{Clarke-CIBC,Clarke-IP} who was the first scholar to undertake the task of developing Whitehead's meretopological ideas. He based his system on the binary relation of \emph{connection}, and the definition of a point different from the original proposal of the English logician. However, as it was later demonstrated by \cite{Biacino-Gerla-CS}, Clarke's contact relation collapses to overlap, and his axioms characterize the atomless complete Boolean algebras. In consequence, Clarke's points as defined in \citep{Clarke-IP} are nothing but ultrafilters. Thus, his approach falls short.

So, let us turn to contact and Boolean contact algebras as a~unifying framework. By a \emph{Boolean contact algebra}\footnote{For an exposition of Boolean contact algebras see \citep{Stell-BCAANATRCC,Bennett-Duntsch-AAT}.} we mean a~Boolean algebra extended with a~binary relation $\con$ of \emph{contact} that satisfies the following constraints:

\begin{gather}
\zero\notcon x\,\label{C0}\tag{C0}\\
x\leq y\wedge x\neq\zero\rarrow x\mathrel{\mathsf{C}} y,\label{C1}\tag{C1}\\
x\mathrel{\mathsf{C}} y\rarrow y\mathrel{\mathsf{C}} x,\label{C2}\tag{C2}\\
x\leq y\rarrow\forall{z\in B}(z\mathrel{\mathsf{C}} x\rarrow z\mathrel{\mathsf{C}} y)\,, \label{C3}\tag{C3}\\
x \con y+z\rarrow x \con y\vee x\con z\,.\label{C4}\tag{C4}
\end{gather}

We extend the inventory of relations by introducing non-tangential inclusion via the expected definition:
\[
x\ll y\iffdef x\notcon-y\,.
\]
The reader will check easily that so defined $\ll$ has properties \eqref{DV1}--\eqref{DV5}. The remaining two De Vries axioms need additional assumptions about~$\con$.

Grzegorczyk's idea to introduce points was somewhat similar to those of Whitehead and De Vries.\footnote{Historically, Grzegorczyk precedes De Vries, yet it is virtually impossible that the two scholars were aware of each other's work.} Take a~region and shrink it till you <<squeeze>> a~point out of it. However, what distinguishes his definition from the other two is that he demanded that every set of regions that is a~candidate for a~point satisfy the following (geometrical in spirit) property: if $x$ and $y$ are regions such that each one overlaps all regions in a~point candidate, then $x$ must touch~$y$ (see Figure~\ref{fig:Grz-representative}). This requirement singles out Grzegorczyk points among De Vries points, as we will see in a~moment.

\begin{figure}[!htbp]
\begin{center}
    \begin{tikzpicture}
    \foreach \y\z in {2.5/10,2.25/20,2/30,1.75/40,1.5/50,1.25/70,1/80,0.75/90,0.5/100}{
    \draw[name path = space,fill=cyan!\z] plot[domain=0:350, smooth cycle] (\x:\y+rnd*0.2);
    }
    \foreach \x in {0,72,...,288}
    {
    \draw [dotted,shorten <= 2pt] (0,0) -- (\x:0.5cm);
    }

    \begin{scope}
    \clip (-2,-4) to [out=90,in=180] ++(2,4) to [out=0,in=90] ++(2,-4) to (-2,-4);
    \fill[yellow,opacity=.7] (-2,-3.5) rectangle (2,0);
    \end{scope}
    \begin{scope}
    \clip (-3,3) to [out=270,in=180] (0,0) to [out=0,in=90] ++(1,3) to cycle;
    \fill[green,opacity=.7] (-3,3) rectangle (1,0);
    \end{scope}
    \end{tikzpicture}\caption{A representative of a~point in the sense of Grzegorczyk: if two regions overlap all elements of the representative, then they must be in contact.}\label{fig:Grz-representative}
     \end{center}
\end{figure}
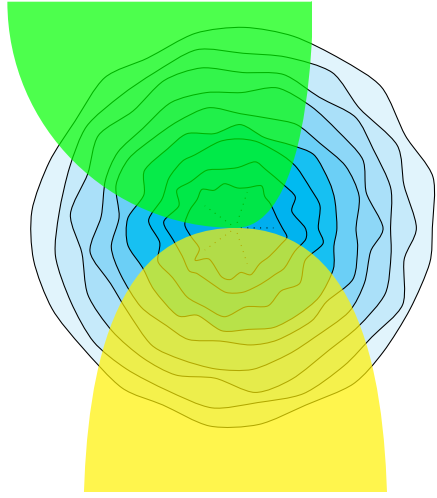

Formally, a~\emph{Grzegorczyk representative of a point} (\emph{G-representative} for short)\footnote{Both the term and its abbreviation adopted from \cite{Biacino-Gerla-CSGWDP}.} in a Boolean contact algebra is a~non-empty set $Q$ of regions such that:
\begin{gather}
\zero\notin Q\,,\tag{r0}\label{r0}\\
\forall{u,v\in Q}(u=v\vee u\ll v \vee v\ll u)\, , \tag{r1}\label{r1}\\
(\forall{u\in Q})(\exists{v\in Q})\; v\ll u\, , \tag{r2}\label{r2}\\
\forall{x,y\in R}\bigl(\forall{u\in Q}(u\overl x\wedge u\overl y) \rarrow x\mathrel{\mathsf{C}} y \bigl)\,, \tag{r3}\label{r3}
\end{gather}
where:
\[
x\overl y\iffdef x\cdot y\neq\zero\,.
\]

It is not hard to see that if we take $\RO(\Real)$ and a~point $r$, then the family $Q\defeq\{(r-\nicefrac{1}{n},r+\nicefrac{1}{n})\mid n\in\Nat^+\}$ is a~G-representative. Of course, different sets may represent the same <<point>>, as for example $Q_e\defeq\{(r-\nicefrac{1}{n},r+\nicefrac{1}{n})\mid n\in\Even^+\}$, $Q_o\defeq\{(r-\nicefrac{1}{n},r+\nicefrac{1}{n})\mid n\in\Odd^+\}$ and $Q$ do ($\Even^+$ and $\Odd^+$ are, respectively, the sets of all positive even integers and of all positive odd integers). It is easy to see that in the case of $r$ (and any other real number) there are uncountably many G-representatives. More generally, if $Q$ is a~G-representative in a~Boolean contact algebra, and $x\in Q$, then the set $\{y\in Q\mid y\leq x\}$ is also a~G-representative. To circumvent the problem of a~unique point identification we declare Grzegorczyk points (G-points) of a Boolean contact algebra~$B$ to be filters generated by G-representatives (whose set is denoted by $\prePt(B)$):
\[
\scrG\in\Grz(B)\iffslim(\exists Q\in\prePt(B))\,\scrG=\{y\in B\mid\exists x\in Q\,y\leq x\}\,.
\]

Let $x\conf\fil$ hold iff region $x$ is in contact with every region in a~filter $\fil$: $(\forall y\in\fil)\,y\con x$. Accordingly, $x\nconf\fil$ iff there is a~region in $\fil$ that is separated from $x$. The reader will easily check that if $\fil$ is round, then:
\begin{equation*}
    x\in\fil\iffslim -x\nconf\fil\,.
\end{equation*}

Interestingly, every Grzegorczyk point is a maximal round filter. Firstly, every G-point~$\scrG$ is a~round filter, for if $x\in\scrG$ and $Q$ generates~$\scrG$, then in $Q$ there is a~region $y\leq x$. But in $Q$ there is $z$ well-inside $y$, so $z$ is well-inside $x$ either. Secondly, in \citep{Gruszczynski-NTP} it was proven that a~round filter $\fil$ satisfies the following condition:
\begin{equation}\tag{$\maltese$}\label{eq:maltese}
(\forall x,y \in B)\,(x\conf\fil\conf y\rarrow x\con y)   \end{equation}
iff the condition \eqref{eq:ends-conditions-2} is also true about~$\fil$. Indeed, if $x$ is separated from $y$, applying \eqref{eq:maltese} in the contraposed form we obtain that either $x\nconf\fil$ or $y\nconf\fil$, and so either $-x$ is an element of $\fil$ or $-y$ is, as required. The reverse implication is proven analogously. As the property \eqref{eq:ends-conditions-2} uniquely identifies maximal round filters, so does~\eqref{eq:maltese}. At the same time, the condition \eqref{r3} for G-representatives, together with the definition of G-points, entail that every Grzegorczyk point must satisfy \eqref{eq:maltese}. So $\Grz(B)\subseteq\End(B)$. Does the other inclusion hold? In general, no. If\label{page:G-not-E} we look back at Figure~\ref{fig:point-at-infinity} we can see a~fragment of a~point at infinity that, in general, does not have to be a~G-point. To see this, imagine that we color the regions of the point with two alternating colors, as in Figure~\ref{fig:G-not-E}.
\begin{figure}[!htbp]
\begin{center}
    \begin{tikzpicture}[scale=0.3]
    \clip (0,-12) rectangle (-12,11.5);
    \foreach \y\z in {9/10,8/20,7/30,6/40,5/50,4/70,3/80,2/90,1/100}{
    \pgfmathparse{isodd(\y)}\ifnum\pgfmathresult=1\def\currcol{cyan}\else\def\currcol{yellow}\fi
    \draw[fill=\currcol!\z] plot[domain=0:350, smooth cycle] (\x:\y+rnd*0.2);
    }
    \foreach \x in {120,180,240}
    {
    \draw [dotted,shorten <= 2pt] (0,0) -- (\x:0.8cm);
    }
    \draw[dotted,thick,gray] (0,-11) rectangle (-12,11.2);
    \end{tikzpicture}
\end{center}\caption{A construction towards showing that not every maximal round filter is a~Grzegorczyk point.}\label{fig:G-not-E}
\end{figure}
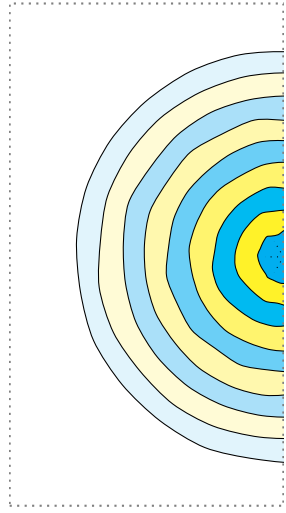
After coloring, we choose only blue stripes, number them with natural numbers, and divide them into two sets: of stripes tagged with even and of stripes tagged with odd numbers, respectively. We can now take the suprema of the first set and the second set to obtain regions that are apart yet overlap every region in the chain we began with. The chain is included in a maximal round filter~$\scrE$, yet $\scrE$ cannot be generated by any G-representative. Precisely because any such a~G-representative would have to be covered by the chain, thus failing to satisfy the condition~\eqref{r3}. We must be careful here as the situation is subtle, so let us repeat: the maximal round filter $\scrE$ must satisfy \eqref{eq:maltese}, since every maximal round filter satisfies the condition; it is only that no G-representative can give rise to~$\scrE$.

Let us show that G-points satisfy the criteria of the method from page \pageref{page:method}. The first step towards demonstrating this was done by Grzegorczyk himself, and more elaborate constructions were delivered by \citet{Gruszczynski-NTP}, and \citet{Gruszczynski-et-al-ASGPFT1,Gruszczynski-et-al-ASGPFT2}. Grzegorczyk demonstrated that his points, together with the Stone-like mapping, form a~topological Hausdorff space, in which his mereology-based separation structures can be represented. Grzegorczyk also maintained that the spaces of his points have the following property: for every point $p$, there exists an infinite strictly decreasing family of open sets such that the intersection of the family is $\{p\}$. Yet this is not true, as there are finite structures that are models of Grzegorczyk axioms, which was proven in the papers by Gruszczyński and Pietruszczak. In those papers, a~class of the so-called \emph{concentric} topological spaces was singled out, which are $T_1$ spaces additionally satisfying the condition \eqref{R1} displayed below on page~\pageref{R1}. Later in \citep{Gruszczynski-et-al-GPFBCA} it was proven that this class forms a~subclass of the so-called \emph{lob-spaces} -- topological spaces with linearly ordered basis at every point (see \citealp{Davis-SLOLB}). The subclass contains only regular spaces; that is, concentric spaces are those lob-spaces that are $T_1$ and regular. Further, both authors proved that every Grzegorczyk structure can be represented as a~subalgebra of the regular open algebra of a~concentric space of Grzegorczyk points. Moreover, it was also proven that there is a~one-to-one correspondence between Grzegorczyk structures that satisfy the countable chain condition and concentric spaces that satisfy the topological version of the condition. As a~result, abstract Grzegorczyk structures obtained concrete representation, and their existence was also established. The latter follows from the fact that, e.g., the real line with the standard Euclidean topology is a~concentric space.

In the BCA setting, a~\emph{Grzegorczyk contact algebra} may be defined as a~Boolean contact algebra that satisfies two additional second-order Grzegorczyk's axioms. The first of them says that every region has a~G-representative (and consequently, a~G-point):
\begin{equation}
    (\forall x\in B)(\exists Q\in\prePt)\,x\in Q\,.\tag{G1}\label{G1}
\end{equation}
According to the second, G-representatives (and so G-points either) exist in those locations of space (understood as the unity of the algebra) where regions touch each other:
\begin{equation}
    x\con y\rarrow(\exists Q\in\prePt)(\forall u\in Q)\,(u\overl x\wedge u\overl y)\,.\tag{G2}\label{G2}
\end{equation}
More precisely, the class of Grzegorczyk contact algebras is determined by axioms \eqref{C0}--\eqref{C3}, \eqref{G1}, \eqref{G2}, as \eqref{C4} is their consequence.

It is provable that the set of all values of the Stone-like mapping $\scrG\colon B\to\power(\Grz B)$ such that $\scrG(x)\defeq\{\scrG\in\Grz(B)\mid x\in\scrG\}$ is a basis, and thus gives rise to a~topological space $\langle\Grz(B),\topo\rangle$. As we wrote above, the key notion to understanding this space is the concept of a~\emph{concentric} space, which is formally defined as a $T_1$ space in which every point $p$ has a local basis $\Basis_p$ of regular open sets such that:
\begin{equation}\tag{R1}\label{R1}
(\forall U,V\in\Basis_p)\,(U=V\vee\Cl U\subseteq V\vee \Cl V\subseteq U)\,.
\end{equation}
The reader will notice that the condition is a~point-based counterpart of \eqref{r1} from page \pageref{r1}.  Every concentric space is a~regular space, yet generally, the converse is not true. For example, the uncountable product of the discrete space $\{0,1\}$ is regular but not concentric \citep{Ruping-AERBNLBTS}.

If $B$ is a~Grzegorczyk contact algebra, then the space $\Grz(B)$ must be a~concentric space. Given any G-point $\scrG$ we know that it has been generated by some G-representative $Q$, and thusly, the family $\Basis_{\scrG}\defeq\{\scrG(x)\mid x\in Q\}$ is a~local basis at the point $\scrG$ that satisfies the condition~\eqref{R1}. The fact that $\Grz(B)$ is $T_1$ is routinely verified, since if $\scrG_1\neq\scrG_2$, then $\scrG_1\nsubseteq\scrG_2$ and $\scrG_2\nsubseteq\scrG_1$ (for G-points are maximal objects). Therefore there is a region $x$ in $\scrG_1$ but not in $\scrG_2$, so $\scrG_1(x)$ is an open set around the point $\scrG_1$ but not around $\scrG_2$.

On the other hand, given a~concentric space $X$, its algebra $\RO(X)$ is a~(complete) Grzegorczyk contact algebra.\footnote{The proof of this fact can be found in \citep{Gruszczynski-NTP} and \citep{Gruszczynski-et-al-ASGPFT2}.}

In \citep{Gruszczynski-NTP} it was shown that every Grzegorczyk contact algebra~$B$ embeds into a~Grzegorczyk contact algebra of a~concentric topological space, and the embedding is an isomorphism in the case of completeness of~$B$.

\begin{figure}[!htbp]
\begin{center}
\begin{tikzcd}[row sep=large,column sep=large]
\GCA\arrow[r, symbol=\ni] &[-3em]B \arrow[r, ->,dashed] \arrow[rd, "\scrG" left=5pt, ->,shift right,hook] & \Grz(B) \arrow[d, dashed] \arrow[r, symbol=\in]& [-3em] \Conc\\
&{}& {\RO(\Grz(B))}\arrow[r, symbol=\in]& [-3em] \GCA^c
\end{tikzcd}\caption{Let $\GCA$ be the class of Grzegorczyk contact algebras, and $\Conc$ the class of concentric topological spaces. Any Grzegorczyk algebra~$B$ embeds densely into the algebra of regular open sets of the concentric space space for $B$.}\label{fig:GCA-to-Conc}
\end{center}
\end{figure}

\begin{figure}[!htbp]
\begin{center}
\begin{tikzcd}[row sep=large,column sep=large]
\GCA^c\arrow[r, symbol=\ni] &[-3em]B \arrow[r, ->,dashed] \arrow[rd, "\scrG" left=5pt, ->,shift right] \arrow[rd, "\scrG^{-1}" right=5pt, <-,shift left] & \Grz(B) \arrow[d, dashed] \arrow[r, symbol=\in]& [-3em] \Conc\\
&{}& {\RO(\Grz(B))}\arrow[r, symbol=\in]& [-3em] \GCA^c
\end{tikzcd}\caption{Any complete Grzegorczyk contact algebra $B$ is indistinguishable from the Grzegorczyk algebra of regular open sets of the concentric space for $B$.}\label{fig:complete-GCA-to-Conc}
\end{center}
\end{figure}

The path from the concentric topological spaces to Grzegorczyk algebras is a bit more complicated, and it was only proven for Grzegorczyk contact algebras and concentric spaces that satisfy, respectively, algebraical and topological versions of the countable chain condition, which has not been circumvented so far. By an antichain of a Boolean algebra we mean, standardly, a~subset of its regions that are pairwise incompatible. In the case of topological spaces, an antichain is a~family of open sets whose intersections are pairwise empty. The countable chain condition is satisfied either by an algebra or a~topological space if any antichain is at most countable.

Firstly, if the condition holds for a~Grzegorczyk algebra~$B$, then its space $\Grz(B)$ satisfies the topological version of the condition, and the algebraical version transfers to $\RO(\Grz(B))$. The first dependence stems from the fact that if every antichain of regions is at most countable and the family of all sets of the form $\scrG(x)$ is a basis for $\Grz(B)$, then the space must also satisfy the condition. If it did not, for an uncountable antichain of its open sets we would find an uncountable antichain of sets of the form $\scrG(x)$, and since:
\[
x\ext y\iffslim\scrG(x)\cap\scrG(y)=\emptyset
\]
the pre-images of $\scrG(x)$s would form an uncountable chain of regions in $B$.

Secondly, it is evident that if a~topological space satisfies the countable chain condition, then its algebra of regular open sets must also satisfy it.

In light of these, it is easily seen that the situations from figures~\ref{fig:GCA-to-Conc} and~\ref{fig:complete-GCA-to-Conc} transfer immediately to those structures that satisfy the condition. Moreover, we can extend the representation to the one from Figure~\ref{fig:Conc-to-GCA}. For complete Grzegorczyk algebras, we have then a~one-to-one correspondence between these that satisfy ccc, and concentric structures that have ccc.

\begin{figure}[!htbp]
\begin{center}
\begin{tikzcd}[row sep=large,column sep=large]
\Conc_{\mathrm{ccc}}\arrow[r, symbol=\ni] &[-3em]X \arrow[r, ->,dashed] \arrow[rd, "h" left=5pt, ->,shift right] \arrow[rd, "h^{-1}" right=5pt, <-,shift left] & \RO(X) \arrow[d, dashed] \arrow[r, symbol=\in]& [-3em] \GCA^c_{\mathrm{ccc}}\\
&{}& {\Grz(\RO(X))}\arrow[r, symbol=\in]& [-3em] \Conc_{\mathrm{ccc}}
\end{tikzcd}\caption{Any concentric space $X$ satisfying the countable chain condition is homeomorphic to the concentric space of the complete Grzegorczyk algebra of regular open sets of $X$ that satisfies the condition either.}\label{fig:Conc-to-GCA}
\end{center}
\end{figure}

To conclude, the results presented let us affirmatively respond to the question: are Grzegorczyk points really points?\footnote{Technical details of all constructions can be found in \citep{Gruszczynski-NTP} and \citep{Gruszczynski-et-al-ASGPFT1,Gruszczynski-et-al-ASGPFT2}.}
As it turns out, thanks to the results for G-points, we can positively answer the main problem of this paper: are there any spaces of \emph{Whitehead points} (in the sense of the method from page~\ref{page:method})?

\section{Spaces of Whitehead points}\label{sec:solution}

In \citep{Biacino-Gerla-CSGWDP} we find the proof that, under some reasonable constraints, the classes of Grzegorczyk points and Whitehead points for a~certain connection structures (a~mereological structures with the contact relation) coincide. In this section, we rephrase the results of \citeauthor{Biacino-Gerla-CSGWDP} in the framework of contact algebras in order to apply their result (together with the results from earlier sections) to the problem of representation theorem for Whitehead points.

As we saw, Grzegorczyk points may be defined as filters, but they can also be characterized as quotients with respect to the covering relation from section~\ref{sec:informal-construction}. In the case\label{page:mutual-cover-Q} of G-representatives we have that if $Q_1$ covers $Q_2$, then $Q_2$ covers $Q_1$. This is a~consequence of two facts: (a) if $Q_1$ does not cover $Q_2$, then there are regions $x\in Q_1$ and $y\in Q_2$ that are separated from each other, and (b) if $Q_1$ covers $Q_2$, then for all $x\in Q_1$ and $y\in Q_2$, $x$ and $y$ are compatible.

Since covering is also transitive and reflexive, it must be an equivalence relation (in the family of G-representatives, but not generally in the family of all abstractive sets), and thus we can say that G-representatives $Q_1$ and $Q_2$ \emph{represent the same location} (in symbols: $Q_1\sim Q_2$) if and only if $Q_1$ covers $Q_2$ (and $Q_2$ covers $Q_1$).

The relation $\sim$ may be recovered from the set of G-points via the following equivalence:
\[
Q_1\sim Q_2\iffslim(\exists\scrG\in\Grz)\,Q_1\cup Q_2\subseteq\scrG\,.
\]

The family of all equivalence classes of the relation $\sim$ in the set of Grzegorczyk representatives:
\[
\Eq\defeq\prePt/_{\mathord{\sim}}\,
\]
may now be treated as the set of points, as there is a~bijective correspondence between elements of $\Eq$ and G-points given by function $f\colon\Eq\to\Grz$ such that $f([Q])\defeq\scrG_Q$, where $\scrG_Q$ is the G-point generated by~$Q$. Thus, Grzegorczyk points can be characterized by the Whiteheadian covering relation.

Let us observe that Whitehead's abstractive sets are sets of regions  that satisfy Grzegorczyk conditions \eqref{r0}, \eqref{r1}, plus non-minimality constraint:
  \begin{equation}\tag{A}\label{A}
    \neg(\exists {x\in A})(\forall {y\in A})\,x\leq y\,.
  \end{equation}
Thus, it is immediate that if the Boolean contact algebra in focus is atomless, then its Grzegorczyk representatives must be abstractive sets. A~less obvious conclusion is that in every atomless contact algebra, every G-representative must be a~Whitehead representative of a point either.
To see this, let us couch---after \citeauthor{Biacino-Gerla-CSGWDP}---a~mathematically satisfactory definition of a~Whitehead representative and a~Whitehead point.

Unlike the covering relation on Grzegorczyk representatives, covering on abstractive sets does not have to be an equivalence relation since it is not generally symmetric. However, it is reflexive and transitive, so it gives rise to the following equivalence relation\footnote{Recall that $A_1\covers A_2$ means $A_1$ \emph{covers} $A_2$.}:
\[
A_1\sim A_2\iffdef A_1\covers A_2\wedge A_2\covers A_1\,.
\]
In the case $A_1\sim A_2$, we say that the objects $A_1$ and $A_2$ are \emph{similar}. The intended meaning of \emph{similarity} is a representation of the same geometrical figure in space. Of course, unlike G-representatives, abstractive sets do not have to represent the same precise \emph{location}, and the idea is to identify those that do. As $\sim$ is an equivalence, we can define---in Whitehead's spirit---geometrical objects as equivalence classes of abstractive sets with respect to similarity, i.e., as elements of the family $\Abs/_{\mathord{\sim}}$.  This family equipped with the following binary relation:
\[
[A_1] \trianglerighteq [A_2]\iffdef A_1\covers A_2
\]
is a~partially ordered set (i.e., $\trianglerighteq$ is reflexive, anti-symmetrical, and transitive).

We can now define Whitehead points and Whitehead representatives. $[A]\in\Abs/_{\mathord{\sim}}$ is a~\emph{Whitehead point} (\emph{W-point}) iff $[A]$ is maximal with respect to~$\trianglerighteq$: for every $[A']\in\Abs/_{\mathord{\sim}}$, if $[A]\trianglerighteq[A']$, then $[A]=[A']$. $A\in\Abs$ is a~\emph{Whitehead representative} of a~point (a~\emph{W-representative}) iff $[A]$ is a~Whitehead point. The set of all Whitehead points and of all Whitehead representatives will be denoted by, respectively, `$\Wthd$' and `$\prePtW$'.

Observe that we can also characterize as W-representatives those abstractive sets that satisfy the following equivalence:
\begin{equation*}
  A\in\prePtW\iffslim(\forall B\in\Abs)\,(A\covers B\rarrow B\covers A)\,.
\end{equation*}
As it was demonstrated in~\citep{Gruszczynski-et-al-GAWPTSC}, the notion of the Whitehead point is consistent, i.e., there are contact algebras with Whitehead points. However, we can still ask: can we prove that there are topological spaces based on Whitehead points obtained in the way described on page~\pageref{page:method}?, and can we find any form of representation theorems for such spaces? Both questions may be answered affirmatively in an indirect way using the result of \citet{Biacino-Gerla-CSGWDP}: under additional assumptions, the set of Whitehead points of a~given contact algebra coincides with the set of Grzegorczyk points.

To prove that every G-point is a~W-point it is enough to show that every G-representative is a~W-representative. This part is relatively easy, and the result from \citep{Biacino-Gerla-CSGWDP} can actually be strengthened to the following (for details, see \citealp{Gruszczynski-et-al-GAWPTSC})
\begin{theorem}\label{th:Q-and-A-is-W}
If $\frB$ is a Boolean contact algebra that satisfies \eqref{DV7} then\/\textup{:} $\frB$ is atomless iff in $\frB$ every G-representative is a W-representative.
\end{theorem}

Proving that every W-representative is a~G-representative is a~bit harder, and the original demonstration of \cite{Biacino-Gerla-CSGWDP} calls for a~small modification. In the class of all abstractive sets of a~given Boolean algebra $B$ we distinguish those that countable abstractive sets can represent. By an $\omega$-abstractive set, we understand an abstractive set $A$ for which there is a~countable abstractive set~$A'$ such that $A$ both covers $A'$ and is covered by~$A'$. Accordingly, $W^\omega$-representatives will be those Whitehead representatives that are $\omega$-abstractive sets. Let $\prePt^\omega_W$ be the set of all $W^\omega$-representatives of a given Boolean contact algebra. We have:
\begin{theorem}
    If $\frB$ is a Boolean contact algebra that satisfies \eqref{DV6} and:
\begin{equation}\tag{C6}\label{C5}
    x\notin\{\zero,\one\}\rarrow x\con-x\,,
\end{equation}
then every $W^\omega$-representative is a~G-representative.
\end{theorem}
The small modification we mentioned is the inclusion of \eqref{C5} in the premises of the theorem. Here, \eqref{C5} is a~region-based version of connectedness, i.e., it says that every non-zero and non-unity region touches its Boolean complement. For a more detailed analysis of this, we again refer the reader to \citep{Gruszczynski-et-al-GAWPTSC}.

In light of the above and the earlier results, we may conclude that the set  $\prePt_G^{\Nat}$ of those Grzegorczyk representatives that countable sets can faithfully represent, we have the equality: $\prePt^{\Nat}_G=\prePt^{\Nat}_W$, and in consequence, $\Grz^\omega=\Wthd^\omega$, where the former set is the set of Grzegorczyk points obtained from the elements of $\prePt^{\Nat}_G$ and the latter the set of Whitehead points obtained from the elements of $\prePt^{\Nat}_W$.

We thus have reached a~point at which we can formulate the following theorem:

\begin{theorem}
Let $B$ be an atomless Boolean contact algebra that satisfies the interpolation axiom \eqref{DV6} and the connectedness axiom \eqref{C5}. Suppose we introduce both definitions of points---by Grzegorczyk and by Whitehead---and extend the axioms with Grzegorczyk postulates \eqref{G1} and \eqref{G2}. Suppose $\Grz^\omega\neq\emptyset$. Let $\langle\Grz,\topo\rangle$ be the concentric topological space for $B$. Then its subspace $\langle\Grz^\omega,\topo^\omega\rangle$ \textup{(}\/where $\topo^\omega\defeq\{\Grz^\omega\cap V\mid V\in\topo\}$\/\textup{)} is a~topological space whose points are W-points.
\end{theorem}


We can also conclude that there are spaces in which both sets of points coincide on the whole space, not only its subspace. To this end, observe that in the case of abstractive sets \emph{covering} is anything but a form of \emph{cofinality} for $\geq$-relation: an abstractive set $A$ covers an abstractive set~$B$ iff $B$ is cofinal with $A$. Putting the dual $\geq$ of \emph{part of} relation in focus, and assuming Axiom of Choice, every chain $C$ in any Boolean contact algebra has a~cofinal well-ordered subchain~$C'$ with respect to $\geq$, where we refer to the dual notion of the well-ordered set by requiring the existence of the maximal element for $\geq$ in every non-empty subset of $C'$. On the other hand, the countable chain condition entails that every infinite well-ordered set of regions must be countable. Therefore:

\begin{theorem}\label{th:2}
Let $B$ be an atomless Boolean contact algebra that satisfies the interpolation axiom \eqref{DV6}, the connectedness axiom \eqref{C5}, and the countable chain condition. Suppose we introduce both definitions of points---by Grzegorczyk and by Whitehead---and we extend the axioms with Grzegorczyk postulates \eqref{G1} and \eqref{G2}.  The concentric topological space $\langle\Grz,\topo\rangle$ for $B$ is a~topological space in which $\Grz=\Grz^\omega$, so it is a~space whose points are W-points.
\end{theorem}

Thanks to the above theorem, we can see that Grzegorczyk and Whitehead points coincide in a large subclass of regular spaces: concentric spaces that satisfy countable chain condition.\footnote{The result concerning the relationship between Grzegorczyk and Whitehead points can be generalized by eliminating the countability assumption. This, however, calls for a stronger, second-order version of \eqref{DV6}. Details, again, can be found in \citep{Gruszczynski-et-al-GAWPTSC}.}

Since the algebra $\RO(\Real^n)$ of regular open subsets of the $n$-dimensional Euclidean space has all the properties from the premises of the theorem above, we can conclude that:
\begin{corollary}
There are spaces of Whitehead points satisfying the requirements of the method from page~\pageref{page:method}.
\end{corollary}

Let us conclude this section with a~strict justification of the difference between De Vries's and Whitehead points mentioned on page \pageref{page:example-DV-not-W}. We know there are structures in which Whitehead points are exactly Grzegorczyk points. Yet on page~\ref{page:G-not-E} we have demonstrated how to construct a~maximally round filter that is not a~G-point. This construction can be carried out in $\Real^2$, which is a~space that satisfies all premises of Theorem~\ref{th:2}. Thus, in $\RO(\Real^2)$ there is a~De Vries point that is not a~Whitehead point.

\section{Summary}

From the intuitions about the perspective space, we have come a long way through the topological representation theorems for Boolean algebras and De Vries algebras, Grzegorczyk contact algebras, to spaces of Whitehead points. Because there are spaces of Grzegorczyk points and Grzegorczyk contact algebras whose G-points coincide with Whitehead points, we concluded that there are topological spaces constructed in the Stone-like manner whose fundamental objects are the English logician's points.

One might say that this is a~roundabout way to show that there are topological spaces built out of Whitehead points. However, to our knowledge, no better way has been found so far. The earlier analyses only presented the way to points via extensive abstraction or compared them to other similar constructions. Yet, none of them pointed out that there are indeed topological spaces of Whitehead points obtained via methods of representation theorems.


The natural questions at this point are: can we generalize the result?, can we drop the reference to Grzegorczyk points and build any representation (or, even better, duality) for Whitehead points directly? With positive answers to these, we may try extending the scrutiny of both Grzegorczyk and Whitehead constructions to algebraic structures weaker than Boolean contact algebras, e.g., (extended) distributive contact lattices \citep{Duntsch-et-al-DCLTR,Ivanova-et-al-DMEDCL}, or Stonian p-ortholattices \citep{Winter-et-al-OTAORS}, to name few.

These, in our opinion, are problems concerning the classical Whitehead construction that has been neglected for too long. The path to understanding what Whitehead points are leads through the realms of logic and mathematics.



\section*{Acknowledgements}
This research was funded by (a) the National Science Center (Poland), grant number~2020/39/B/HS1/00216.

For the purpose of Open Access, the author has applied a CC-BY public copyright license to any Author Accepted Manuscript (AAM) version arising from this submission.

\bibliographystyle{apalike}
\providecommand{\noop}[1]{}

\end{document}

However, we can make up for these two in a~relatively easy way.

After~\cite{Shchepin-RFNNS}, let us call a~topological space \emph{$\kappa$-normal}\footnote{Another term, use by \citealp{Bennett-Duntsch-AAT}, is `weakly-normal'.} if and only if any two disjoint non-empty regular closed sets can be separated by (regular) open sets. This is equivalent to the following:
\[
(\forall U,V\in\RO(X)^+)\,(\Cl U\subseteq V\rarrow(\exists M\in\RO(X))\,\Cl U\subseteq M\wedge\Cl M\subseteq V)\,.
\]
The interpolation axiom \eqref{DV6} is a~counterpart of $\kappa$-normality under the standard topological interpretation of the contact relation \citep{Duntsch-et-al-RTBCA}. Thanks to \cite{Vakarelov-et-al-ANOPSACBM} we know that every Boolean contact algebra that satisfies the interpolation axiom is isomorphic to a~dense substructure of the regular open algebra of a~completely regular\footnote{A space $X$ is \emph{completely regular} if any closed subset $A$ of $X$ and a~point $p$ beyond $A$ can be separated by a~continuous function $f\colon X\to[0,1]$.} $T_1$ space with the standard topological interpretation of contact. The result does not transform immediately to Grzegorczyk contact alegbras, since \citeauthor{Vakarelov-et-al-ANOPSACBM} use different objects as points, the so called~\emph{clusters}.